\documentclass[12pt,reqno]{amsart}
\usepackage{dsfont,amsmath,amsfonts,amscd,amssymb,graphicx,mathrsfs,eufrak,enumitem,euscript,colonequals,upgreek,pifont}
\usepackage[all,arc,curve,color,frame]{xy}
\usepackage[usenames]{color}

\usepackage{setspace}
\usepackage{array,multirow,booktabs,longtable}

\setlist[enumerate]{format=\normalfont}

\usepackage{pifont}

\newcommand\Curve{\mathrm{C}}

\newcommand{\tocspace}{0.1ex}
\let\oldtocsection=\tocsection
\let\oldtocsubsection=\tocsubsection
\let\oldtocsubsubsection=\tocsubsubsection
\renewcommand{\tocsection}[3]{\hspace{0em}\oldtocsection{#1}{#2}{#3}\vspace{\tocspace}}
\renewcommand{\tocsubsection}[3]{ \hspace{1em} \oldtocsubsection{#1}{\small{#2}}{\small{#3}}\vspace{\tocspace} }
\renewcommand{\tocsubsubsection}[3]{\hspace{2em}\oldtocsubsubsection{#1}{\small{#2}}{\small{#3}}}

\newcommand{\marginparstretch}{0.6}
\let\oldmarginpar\marginpar
\renewcommand\marginpar[1]{\-\oldmarginpar[\framebox{\setstretch{\marginparstretch}\begin{minipage}{\marginparwidth}{\raggedleft\tiny #1}\end{minipage}}]{\framebox{\setstretch{\marginparstretch}\begin{minipage}{\marginparwidth}{\raggedright\tiny #1}\end{minipage}}}}

\usepackage[colorlinks]{hyperref}
\usepackage{tikz,mathrsfs}
\usetikzlibrary{arrows,decorations.pathmorphing,decorations.pathreplacing,positioning,shapes.geometric,shapes.misc,decorations.markings,decorations.fractals,calc,patterns}

\tikzset{
        cvertex/.style={circle,draw=black,inner sep=1pt,outer sep=3pt},
        vertex/.style={circle,fill=black,inner sep=1pt,outer sep=3pt},
        star/.style={circle,fill=yellow,inner sep=0.75pt,outer sep=0.75pt},
        tvertex/.style={inner sep=1pt,font=\scriptsize},
        gap/.style={inner sep=0.5pt,fill=white}}

\tikzset{
        DB/.style={circle,draw=black,circle,fill=white,inner sep=0pt, minimum size=5pt},
        DW/.style={circle,draw=black,fill=black,inner sep=0pt, minimum size=5pt},
        dotted/.style={circle,draw=black,densely dotted,fill=white,inner sep=0pt, minimum size=5pt},
}

\tikzstyle{mybox} = [draw=black, fill=blue!10, very thick,
    rectangle, rounded corners, inner sep=10pt, inner ysep=20pt]
\tikzstyle{boxtitle} =[fill=blue!50, text=white,rectangle,rounded corners]

\newtheorem{thm}{Theorem}[section]
\newtheorem{prop}[thm]{Proposition}

\newtheorem{defin}[thm]{Definition}

\theoremstyle{definition} 

\newtheorem{example}[thm]{Example}

\newtheorem{remark}[thm]{Remark}

\numberwithin{equation}{section}

\newcommand\citetype[1]{}

\newcommand{\m}{\mathfrak{m}}

\def\coh{\mathop{\rm coh}\nolimits}

\def\RHom{\mathop{\rm {\bf R}Hom}\nolimits}
\def\End{\mathop{\rm End}\nolimits}
\def\Ext{\mathop{\rm Ext}\nolimits}

\def\Spec{\mathop{\rm Spec}\nolimits}

\def\Auteq{\mathop{\rm{Auteq}}\nolimits}

\def\Db{\mathop{\rm{D}^b}\nolimits}

\def\flop{{\sf{F}}}

\def\Id{\mathop{\rm{Id}}\nolimits}

\newcommand{\con}{\mathrm{con}}

\newcommand{\CA}{\mathrm{A}_{\con}}

\def\RHom{{\rm{\bf R}Hom}}

\newcommand\twistGen{{\sf Twist}}

\newcommand{\scrA}{\EuScript{A}}

\newcommand{\scrE}{\EuScript{E}}

\newcommand{\scrO}{\EuScript{O}}

\newcommand{\scrS}{\EuScript{S}}

\newcommand{\scrV}{\EuScript{V}}
\newcommand{\scrZ}{\EuScript{Z}}

\newcommand{\sigE}{{^{\upsigma}\kern -1pt\scrE}}
\newcommand{\nsigE}{{^{{\scriptsize-}\upsigma}\kern -1pt\scrE}}

\def\Id{\mathop{\rm{Id}}\nolimits}

\def\Db{\mathop{\rm{D}^b}\nolimits}

\numberwithin{equation}{section}

\makeatletter
\newcommand*\bigcdot{\mathpalette\bigcdot@{.5}}
\newcommand*\bigcdot@[2]{\mathbin{\vcenter{\hbox{\scalebox{#2}{$\m@th#1\bullet$}}}}}
\makeatother

\begin{document}

\title[Autoequivalences for $3$-fold Flops]{\textsc{Autoequivalences for 3-fold Flops:\ an overview}}
\author{Michael Wemyss}
\address{Michael Wemyss, The Mathematics and Statistics Building, University of Glasgow, University Place, Glasgow, G12 8QQ, UK.}
\email{michael.wemyss@glasgow.ac.uk}
\begin{abstract}
This is an overview article, based on my 2018 Kinosaki lecture, that surveys and announces work on 3-fold flopping contractions, their affine combinatorics, stability conditions, tilting bundles and autoequivalences. Some first applications are given.
\end{abstract}
\thanks{The author was supported by EPSRC grants~EP/R009325/1 and EP/R034826/1.}
\maketitle
\parindent 20pt
\parskip 0pt

The setting of the lecture are $3$-fold flopping contractions $X\to \Spec R$, where $X$ has at worst Gorenstein terminal singularities.  Many of the statements globalise to projective varieties.

\section{Combinatorics}\label{L0}

Dynkin combinatorics, finite and affine, are well known.  It turns out that the combinatorics of flops does not quite fit into this classical picture: to describe the birational geometry needs new structures and phenomena.

\subsection{Single Vertex Example}\label{S1.1}  This section constructs the 2-sphere $S^2$, minus $6$ points.  This is, of course, not the most difficult mathematical object to produce. The point here is that the construction goes via Dynkin diagrams.

Just the following two facts are required.
\begin{enumerate}
\item ADE Dynkin diagrams are the following graphs, where the vertices have been labelled with the data of  the rank of the highest root; geometrically, this is the fundamental cycle in Kleinian singularities.
\[
\begin{tikzpicture}
\node at (0,0.25) {$A_n$};
\node at (0,-0.25) {$\scriptstyle n\geq 1$};
\node at (5.01,0) {$
\begin{tikzpicture}
 \node (0) at (0,0) [DB] {};
 \node (1) at (1,0) [DB] {};
 \node (2) at (2,0) [DB] {};
 \node (3) at (4,0) [DB] {};
 \node (4) at (5,0)[DB] {};
 \node at (3,0) {$\cdots$};
 \node (0a) at (0,-0.3) {$\scriptstyle 1$};
 \node (1a) at (1,-0.3) {$\scriptstyle 1$};
 \node (2a) at (2,-0.3) {$\scriptstyle 1$};
 \node (3a) at (4,-0.3) {$\scriptstyle 1$};
 \node (4a) at (5,-0.3) {$\scriptstyle 1$};
\draw [-] (0) -- (1);
\draw [-] (1) -- (2);
\draw [-] (2) -- (2.6,0);
\draw [-] (3.4,0) -- (3);
\draw [-] (3) -- (4);
\end{tikzpicture}$};
\node at (0,-1.25) {$D_n$};
\node at (0,-1.75) {$\scriptstyle n\geq 4$};
\node at (5.01,-1.5) {\begin{tikzpicture}
 \node (0) at (0,0) [DB] {};
 \node (1) at (1,0) [DB] {};
 \node (1b) at (1,0.75) [DB] {};
 \node (2) at (2,0) [DB] {};
 \node (3) at (4,0) [DB] {};
 \node (4) at (5,0)[DB] {};
 \node at (3,0) {$\cdots$};
 \node (0a) at (0,-0.3) {$\scriptstyle 1$};
 \node (1a) at (1,-0.3) {$\scriptstyle 2$};
 \node (1ba) at (0.75,0.75) {$\scriptstyle 1$};
 \node (2a) at (2,-0.3) {$\scriptstyle 2$};
 \node (3a) at (4,-0.3) {$\scriptstyle 2$};
 \node (4a) at (5,-0.3) {$\scriptstyle 1$};
\draw [-] (0) -- (1);
\draw [-] (1) -- (2);
\draw [-] (2) -- (2.6,0);
\draw [-] (3.4,0) -- (3);
\draw [-] (3) -- (4);
\draw [-] (1) -- (1b);
\end{tikzpicture}};
\node at (0,-3) {$E_6$};
\node at (4.1,-3) {$\begin{tikzpicture}[scale=0.8]
 \node (m1) at (-1,0) [DB] {};
 \node (0) at (0,0) [DB] {};
 \node (1) at (1,0) [DB] {};
 \node (1b) at (1,0.75) [DB] {};
 \node (2) at (2,0) [DB] {};
 \node (3) at (3,0) [DB] {};
\node (m1a) at (-1,-0.3) {$\scriptstyle 1$};
 \node (0a) at (0,-0.3) {$\scriptstyle 2$};
 \node (1a) at (1,-0.3) {$\scriptstyle 3$};
 \node (1ba) at (0.75,0.75) {$\scriptstyle 2$};
 \node (2a) at (2,-0.3) {$\scriptstyle 2$};
 \node (3a) at (3,-0.3) {$\scriptstyle 1$};
\draw [-] (m1) -- (0);
\draw [-] (0) -- (1);
\draw [-] (1) -- (2);
\draw [-] (2) -- (3);
\draw [-] (1) -- (1b);
\end{tikzpicture}$};
\node at (0,-4.5) {$E_7$};
\node at (4.5,-4.5) {$\begin{tikzpicture}[scale=0.8]
 \node (m1) at (-1,0) [DB] {};
 \node (0) at (0,0) [DB] {};
 \node (1) at (1,0) [DB] {};
 \node (1b) at (1,0.75) [DB] {};
 \node (2) at (2,0) [DB] {};
 \node (3) at (3,0) [DB] {};
 \node (4) at (4,0) [DB] {};
\node (m1a) at (-1,-0.3) {$\scriptstyle 2$};
 \node (0a) at (0,-0.3) {$\scriptstyle 3$};
 \node (1a) at (1,-0.3) {$\scriptstyle 4$};
 \node (1ba) at (0.75,0.75) {$\scriptstyle 2$};
 \node (2a) at (2,-0.3) {$\scriptstyle 3$};
 \node (3a) at (3,-0.3) {$\scriptstyle 2$};
 \node (4a) at (4,-0.3) {$\scriptstyle 1$};
\draw [-] (m1) -- (0);
\draw [-] (0) -- (1);
\draw [-] (1) -- (2);
\draw [-] (2) -- (3);
\draw [-] (3) -- (4);
\draw [-] (1) -- (1b);
\end{tikzpicture}$};
\node at (0,-6) {$E_8$};
\node at (4.9,-6) {$\begin{tikzpicture}[scale=0.8]
 \node (m1) at (-1,0) [DB] {};
 \node (0) at (0,0) [DB] {};
 \node (1) at (1,0) [DB] {};
 \node (1b) at (1,0.75) [DB] {};
 \node (2) at (2,0) [DB] {};
 \node (3) at (3,0) [DB] {};
 \node (4) at (4,0) [DB] {};
  \node (5) at (5,0) [DB] {};
\node (m1a) at (-1,-0.3) {$\scriptstyle 2$};
 \node (0a) at (0,-0.3) {$\scriptstyle 4$};
 \node (1a) at (1,-0.3) {$\scriptstyle 6$};
 \node (1ba) at (0.75,0.75) {$\scriptstyle 3$};
 \node (2a) at (2,-0.3) {$\scriptstyle 5$};
 \node (3a) at (3,-0.3) {$\scriptstyle 4$};
 \node (4a) at (4,-0.3) {$\scriptstyle 3$};
 \node (5a) at (5,-0.3) {$\scriptstyle 2$};
\draw [-] (m1) -- (0);
\draw [-] (0) -- (1);
\draw [-] (1) -- (2);
\draw [-] (2) -- (3);
\draw [-] (3) -- (4);
\draw [-] (4) -- (5);
\draw [-] (1) -- (1b);
\end{tikzpicture}$};
\end{tikzpicture}
\]
\item  ADE Dynkin diagrams carry a canonical involution $\upiota$, sometimes called the \emph{Dynkin involution}.  In Type $A_n$ it acts as a reflection in the centre point of the chain (which may or may not be a vertex):
\[
\begin{tikzpicture}[>=stealth]
 \node (1) at (1,0) [DB] {};
 \node (2) at (2,0) [DB] {};
 \node (3) at (4,0) [DB] {};
 \node (4) at (5,0)[DB] {};
\draw [-] (1) -- (2);
\draw [-] (2) -- (2.6,0);
\draw [-] (3.4,0) -- (3);
\draw [-] (3) -- (4);
\draw[densely dotted] (3,0.5)--(3,-0.5);
\draw[<->, bend left,densely dotted] (2.75,0.5) to (3.25,0.5);
\end{tikzpicture}
\]
When $n\geq 4$ and $n$ is odd, the involution $\upiota$ acts on $D_{n}$ by permuting the left hand branches, and fixing all other vertices:
\[
\begin{tikzpicture}[>=stealth]
 \node (0) at (0,0) [DB] {};
 \node (1) at (1,0) [DB] {};
 \node (1b) at (1,0.75) [DB] {};
 \node (2) at (2,0) [DB] {};
 \node (3) at (4,0) [DB] {};
 \node (4) at (5,0)[DB] {};
 \node at (3,0) {$\cdots$};
\draw [-] (0) -- (1);
\draw [-] (1) -- (2);
\draw [-] (2) -- (2.6,0);
\draw [-] (3.4,0) -- (3);
\draw [-] (3) -- (4);
\draw [-] (1) -- (1b);
\draw[<->, bend left,densely dotted] (0) to (1b);
\end{tikzpicture}
\]
In type $E_6$ the involution $\upiota$ acts as a reflection:
\[
\begin{tikzpicture}[>=stealth]
 \node (m1) at (-1,0) [DB] {};
 \node (0) at (0,0) [DB] {};
 \node (1) at (1,0) [DB] {};
 \node (1b) at (1,0.75) [DB] {};
 \node (2) at (2,0) [DB] {};
 \node (3) at (3,0) [DB] {};
\draw [-] (m1) -- (0);
\draw [-] (0) -- (1);
\draw [-] (1) -- (2);
\draw [-] (2) -- (3);
\draw [-] (1) -- (1b);
\draw[densely dotted] (1,-0.5)--(1,1.25);
\draw[<->, bend left,densely dotted] (0.75,1) to (1.25,1);
\end{tikzpicture}
\]
In all other cases, $\upiota$ is the identity.
\end{enumerate}

I claim, given any ADE Dynkin digram, and any choice of vertex, that it is possible to construct a sphere, minus a certain number of points. The dependence between the number of points and the initial choice of vertex is explained later, in Remarks \ref{which 3 not relevant} and \ref{which ell not relevant}.

\begin{example}\label{E6 3 example}
I demonstrate how to construct $S^2\backslash \{ 6 \mathrm{pts}\}$, from the initial data of the $E_6$ Dynkin diagram, and choice of middle vertex, labelled $3$.
\[
\begin{tikzpicture}[>=stealth]
 \node (m1) at (-1,0) [DB] {};
 \node (0) at (0,0) [DB] {};
 \node (1) at (1,0) [DW] {};
 \node (1b) at (1,0.75) [DB] {};
 \node (2) at (2,0) [DB] {};
 \node (3) at (3,0) [DB] {};
\draw [-] (m1) -- (0);
\draw [-] (0) -- (1);
\draw [-] (1) -- (2);
\draw [-] (2) -- (3);
\draw [-] (1) -- (1b);
\end{tikzpicture}
\]
\textbf{Step 1.} First add the extended vertex, which is labelled $1$, and also shaded, to obtain
\newcommand{\EbaseA}[1][]{%
\begin{tikzpicture}[scale=0.21]
\node at (0,0) [DB] {};
\node at (1,0) [DB] {};
\node at (2,0) [DW] {};
\node at (2,1) [DB] {};
\node at (2,2) [DW] {};
\node at (3,0) [DB] {};
\node at (4,0) [DB] {};
\end{tikzpicture}
}
\newcommand{\EbaseB}[1][]{%
\begin{tikzpicture}[scale=0.21]
\node at (0,0) [DB] {};
\node at (1,0) [DB] {};
\node at (2,0) [DW] {};
\node at (2,1) [DW] {};
\node at (2,2) [DB] {};
\node at (3,0) [DB] {};
\node at (4,0) [DB] {};
\end{tikzpicture}
}
\[
\EbaseA
\]
\textbf{Step 2.}  Iterate Dynkin involutions.  Choose one of the shaded vertices, say the one labelled $1$.  We then temporarily delete this vertex, and apply the Dynkin involution to the remaining $E_6$ arrangement
\[
\begin{tikzpicture}[scale=0.21]
\node at (0,0) [DB] {};
\node at (1,0) [DB] {};
\node at (2,0) [DW] {};
\node at (2,1) [DB] {};
\node at (2,2) [dotted] {};
\node at (3,0) [DB] {};
\node at (4,0) [DB] {};
\end{tikzpicture}
\]
The $E_6$ involution fixes the shaded vertex, and so nothing happens.  Since this `move' involved choosing the shaded vertex labelled $1$, we record this as
\[
\begin{tikzpicture}
\draw[densely dotted,->] (-4.5,0) -- (0.5,0);
{\foreach \i in {1}
\filldraw[fill=white,draw=black] (-2*\i,0) circle (2pt);
}
{\foreach \i in {2,1}
\node at (-2*\i+1,0.5) {\EbaseA};
}
\node at (-2,-0.25) {$\scriptstyle 1$};
\end{tikzpicture}
\]
Next, choose the other shaded vertex, the one labelled $3$.  We then temporarily delete this vertex, and apply the Dynkin involution to the remaining vertices, which form a disjoint union of Type $A$ arrangements
\[
\begin{tikzpicture}[scale=0.21,>=stealth]
\node (a) at (0,0) [DB] {};
\node (b) at (1,0) [DB] {};
\node (c) at (2,0) [dotted] {};
\node (d)  at (2,1) [DB] {};
\node (e) at (2,2) [DW] {};
\node (f) at (3,0) [DB] {};
\node (g) at (4,0) [DB] {};
\draw[<->] (-0.15,-0.7) to (1.15,-0.7);
\draw[<->] (2.85,-0.7) to (4.15,-0.7);
\draw[<->] (2.7,0.85) to (2.7,2.15);
\end{tikzpicture}
\]
This moves the shaded vertex down.  As this `move' involved choosing the shaded vertex labelled $3$, we record this move by extending the previous picture to right:
\[
\begin{tikzpicture}
\draw[densely dotted,->] (-4.5,0) -- (2.5,0);
{\foreach \i in {1,0}
\filldraw[fill=white,draw=black] (-2*\i,0) circle (2pt);
}
{\foreach \i in {2,1}
\node at (-2*\i+1,0.5) {\EbaseA};
}
{\foreach \i in {0}
\node at (-2*\i+1,0.5) {\EbaseB};
}
\node at (-2,-0.25) {$\scriptstyle 1$};
\node at (0,-0.25) {$\scriptstyle 3$};
\end{tikzpicture}
\]
Next, choose the \emph{other} vertex, the one now labelled $2$.  Temporarily deleting it, applying the Dynkin involution, then nothing happens.  We again extend the picture to the right.  Continuing in this way, we arrive at:
\[
\begin{tikzpicture}
\draw[densely dotted,->] (-4.5,0) -- (6.5,0);
\node at (6.75,0) {$\mathbb{R}$};
{\foreach \i in {-3,-2,-1,0,1,2}
\filldraw[fill=white,draw=black] (-2*\i,0) circle (2pt);
}
{\foreach \i in {2,1}
\node at (-2*\i+1,0.5) {\EbaseA};
}
{\foreach \i in {0,-1}
\node at (-2*\i+1,0.5) {\EbaseB};
}
\node at (5,0.5) {\EbaseA};
\node at (-4,-0.25) {$\scriptstyle 3$};
\node at (-2,-0.25) {$\scriptstyle 1$};
\node at (0,-0.25) {$\scriptstyle 3$};
\node at (2,-0.25) {$\scriptstyle 2$};
\node at (4,-0.25) {$\scriptstyle 3$};
\node at (6,-0.25) {$\scriptstyle 1$};
\end{tikzpicture}
\]
and observe that this repeats again and again.\\
\textbf{Step 3.} Complexify this hyperplane arrangement, to obtain
\[
\begin{tikzpicture}
\draw[densely dotted,->] (-4.5,0) -- (6.5,0);
\node at (6.75,0) {$\mathbb{R}$};
\draw[densely dotted,->] (1,-0.5) -- (1,1);
\node at (1.4,1) {$\mathrm{i}\mathbb{R}$};
{\foreach \i in {-3,-2,-1,0,1,2}
\filldraw[fill=white,draw=black] (-2*\i,0) circle (2pt);
}
\node at (-4,-0.25) {$\scriptstyle 3$};
\node at (-2,-0.25) {$\scriptstyle 1$};
\node at (0,-0.25) {$\scriptstyle 3$};
\node at (2,-0.25) {$\scriptstyle 2$};
\node at (4,-0.25) {$\scriptstyle 3$};
\node at (6,-0.25) {$\scriptstyle 1$};
\end{tikzpicture}
\]
\textbf{Step 4.} The `class group' $\mathbb{Z}$ acts on the complexification.  This is purely combinatorial: the class group likes the number one, and its generator acts by moving any wall labelled $1$ to the next wall labelled one.  That is to say, $\mathbb{Z}$ acts via translation
\[
\begin{tikzpicture}
\filldraw[gray!10!white] (-2.5,-0.5) -- (-2.5,1) -- (5.5,1)--(5.5,-0.5) --cycle;
\draw[thick,->] (-2.5,0.5) -- (5.5,0.5);
\draw (-2.5,1) -- (-2.5,-0.5);
\draw[dotted] (5.5,-0.5) -- (5.5,1);
\draw[densely dotted,->] (1,-0.5) -- (1,1);
\draw[densely dotted,->] (-4.5,0) -- (6.5,0);
{\foreach \i in {-3,-2,-1,0,1,2}
\filldraw[fill=white,draw=black] (-2*\i,0) circle (2pt);
}
\node at (-4,-0.25) {$\scriptstyle 3$};
\node at (-2,-0.25) {$\scriptstyle 1$};
\node at (0,-0.25) {$\scriptstyle 3$};
\node at (2,-0.25) {$\scriptstyle 2$};
\node at (4,-0.25) {$\scriptstyle 3$};
\node at (6,-0.25) {$\scriptstyle 1$};
\end{tikzpicture}
\]
where for clarity we have shaded the fundamental region of the action. \\
\textbf{Step 5.} Consider the quotient $(\mathbb{C}\backslash \{\rm{pts}\})/\mathbb{Z}$.  This just means we identify the rightmost and leftmost edges of the fundamental region in the above picture, and obtain
\[
\begin{array}{c}
\begin{tikzpicture}[scale=0.85]
\draw (1,-1) -- (1,1);
\draw (-1,-1) -- (-1,1);
\draw (0,1) ellipse (1cm and 0.2cm);
\draw[densely dotted] (1,-1) arc (0:180:1cm and 0.2cm);
\draw (1,-1) arc (0:-180:1cm and 0.2cm);
\draw[gray,densely dotted] (1,0) arc (0:180:1cm and 0.2cm);
\draw[densely dotted] (1,0) arc (0:-180:1cm and 0.2cm)
coordinate[pos=0.71] (A) coordinate[pos=0.57] (B) coordinate[pos=0.43] (C) coordinate[pos=0.29] (D);
\draw[draw=none] (1,-0.2) arc (0:-180:1cm and 0.2cm)
coordinate[pos=0.71] (a) coordinate[pos=0.57] (b) coordinate[pos=0.43] (c) coordinate[pos=0.29] (d);
\filldraw[fill=white,draw=black] (A) circle (1.5pt);
\filldraw[fill=white,draw=black] (B) circle (1.5pt);
\filldraw[fill=white,draw=black] (C) circle (1.5pt);
\filldraw[fill=white,draw=black] (D) circle (1.5pt);
\node at (a) {$\scriptstyle 1$};
\node at (b) {$\scriptstyle 3$};
\node at (c) {$\scriptstyle 2$};
\node at (d) {$\scriptstyle 3$};
\end{tikzpicture}
\end{array}
\]
Topologically, this is just the two-sphere, minus 6 points:
\[
\begin{array}{c}
\begin{tikzpicture}
\draw[gray,densely dotted] (1,0) arc (0:180:1cm and 0.2cm);
\draw[densely dotted] (1,0) arc (0:-180:1cm and 0.2cm)
coordinate[pos=0.71] (A) coordinate[pos=0.57] (B) coordinate[pos=0.43] (C) coordinate[pos=0.29] (D);
\draw ([shift=(-84:1cm)]0,0) arc (-84:84:1cm)  
[bend left] to (96:1cm)
arc (96:264:1cm)
[bend left] to cycle; 
\draw[densely dotted] ([shift=(86.75:1cm)]0,0) arc  (86.75:94:1cm);
\draw[densely dotted] ([shift=(-86.75:1cm)]0,0) arc  (-86.75:-94:1cm);
\draw[draw=none] (1,-0.2) arc (0:-180:1cm and 0.2cm)
coordinate[pos=0.71] (a) coordinate[pos=0.57] (b) coordinate[pos=0.43] (c) coordinate[pos=0.29] (d);
\filldraw[fill=white,draw=black] (A) circle (1.5pt);
\filldraw[fill=white,draw=black] (B) circle (1.5pt);
\filldraw[fill=white,draw=black] (C) circle (1.5pt);
\filldraw[fill=white,draw=black] (D) circle (1.5pt);
\node at (a) {$\scriptstyle 1$};
\node at (b) {$\scriptstyle 3$};
\node at (c) {$\scriptstyle 2$};
\node at (d) {$\scriptstyle 3$};
\end{tikzpicture}
\end{array}
\]
\end{example}
The above construction is somewhat counter-intuitive.  I know of no other way of immediately obtaining the number four (the number of holes on the equator) from the choice of the vertex $3$ in the $E_6$ Dynkin diagram.  A birational-geometric interpretation of the holes on the equator, and their labels, will be given later.

\begin{remark}\label{which 3 not relevant}
Say we perform the same calculation as in Example~\ref{E6 3 example} using instead the initial input of: 
\begin{enumerate}
\item The $E_7$ Dynkin diagram.
\item A shaded vertex labelled $3$. 
\end{enumerate}
It turns out that we get exactly the same answer as before, namely $S^2\backslash \{ 6 \mathrm{pts}\}$.  Even more remarkably, the same labels $1,3,2,3$ again appear on the equator.  Indeed, the calculation is now:
\newcommand{\EsevenbaseA}[1][]{%
\begin{tikzpicture}[scale=0.21]
\node at (-1,0) [DW] {};
\node at (0,0) [DB] {};
\node at (1,0) [DB] {};
\node at (2,0) [DB] {};
\node at (2,1) [DB] {};
\node at (3,0) [DW] {};
\node at (4,0) [DB] {};
\node at (5,0) [DB] {};
\end{tikzpicture}
}
\newcommand{\EsevenbaseB}[1][]{%
\begin{tikzpicture}[scale=0.21]
\node at (-1,0) [DB] {};
\node at (0,0) [DB] {};
\node at (1,0) [DB] {};
\node at (2,0) [DB] {};
\node at (2,1) [DW] {};
\node at (3,0) [DW] {};
\node at (4,0) [DB] {};
\node at (5,0) [DB] {};
\end{tikzpicture}
}
\newcommand{\EsevenbaseC}[1][]{%
\begin{tikzpicture}[scale=0.21]
\node at (-1,0) [DB] {};
\node at (0,0) [DB] {};
\node at (1,0) [DW] {};
\node at (2,0) [DB] {};
\node at (2,1) [DW] {};
\node at (3,0) [DB] {};
\node at (4,0) [DB] {};
\node at (5,0) [DB] {};
\end{tikzpicture}
}
\newcommand{\EsevenbaseD}[1][]{%
\begin{tikzpicture}[scale=0.21]
\node at (-1,0) [DB] {};
\node at (0,0) [DB] {};
\node at (1,0) [DW] {};
\node at (2,0) [DB] {};
\node at (2,1) [DB] {};
\node at (3,0) [DB] {};
\node at (4,0) [DB] {};
\node at (5,0) [DW] {};
\end{tikzpicture}
}
\[
\begin{tikzpicture}
\filldraw[gray!10!white] (-2.5,-0.5) -- (-2.5,1) -- (5.5,1)--(5.5,-0.5) --cycle;
\draw[densely dotted,->] (-4.5,0) -- (6.5,0);
\node at (6.75,0) {$\mathbb{R}$};
{\foreach \i in {-3,-2,-1,0,1,2}
\filldraw[fill=white,draw=black] (-2*\i,0) circle (2pt);
}

\node at (-3,0.5) {\EsevenbaseA};
\node at (-1,0.5) {\EsevenbaseA};
\node at (1,0.5) {\EsevenbaseB};
\node at (3,0.5) {\EsevenbaseC};
\node at (5,0.5) {\EsevenbaseD};
\node at (-4,-0.25) {$\scriptstyle 3$};
\node at (-2,-0.25) {$\scriptstyle 1$};
\node at (0,-0.25) {$\scriptstyle 3$};
\node at (2,-0.25) {$\scriptstyle 2$};
\node at (4,-0.25) {$\scriptstyle 3$};
\node at (6,-0.25) {$\scriptstyle 1$};
\end{tikzpicture}
\]
Ditto for $E_8$.  Both choices of vertex labelled by $3$ give, as before,  $S^2\backslash \{ 6 \mathrm{pts}\}$ with the labels $1,3,2,3$ on the equator.  

\textbf{Upshot:} the labelled topological output is \emph{a property of the number three}, and does not rely on the ambient Dynkin diagram.
\end{remark}

\begin{remark}\label{which ell not relevant}
The above Remark also holds generally.  Pick a vertex in any ADE Dynkin diagram, which is necessarily labelled $\ell$ for some $1\leq \ell\leq 6$.  Again the outcome of the calculation in Example~\ref{E6 3 example} for this choice turns out to only depend on $\ell$, and not on the ambient ADE Dynkin diagram.

The following table summarises, for choice of vertex labelled $\ell$, the number of holes $N$ that appear on the equator, so that the topological space constructed is $S^2\backslash\{N+2\}$.  It also summarises the labels of those holes on the equator.
\[
\begin{tabular}{ccccl}
\toprule
$\ell$&&$N$&&\textnormal{Labels on Equator}\\
\midrule
$1$ && $1$ && $1$ \\
$2$ && $2$ && $1$,$2$ \\
$3$ && $4$ && $1$,$3$,$2$,$3$ \\
$4$ && $6$ && $1$,$4$,$3$,$2$,$3$,$4$ \\
$5$ && $10$ && $1$,$5$,$4$,$3$,$5$,$2$,$5$,$3$,$4$,$5$ \\
$6$ && $12$ && $1$,$6$,$5$,$4$,$3$,$5$,$2$,$5$,$3$,$4$,$5$,$6$ \\
\bottomrule
\end{tabular}
\]
The jump in the number $N$ between $\ell=4$ and $\ell=5$ is unexpected, as is the emergence of the number $5$ in the middle of the labels (for $\ell=5,6$) in what is otherwise a simple pattern.
\end{remark}

\subsection{Two Vertex Example}
Given the choice of a single node in an ADE Dynkin diagram, the previous section first constructed an infinite hyperplane arrangement in $\mathbb{R}^1$, then complexified, then took the quotient by a naturally defined $\mathbb{Z}$-action.

This generalises.  For example, given the choice of two nodes in any ADE Dynkin diagram, using a similar construction it is possible to first produce an infinite hyperplane arrangement in $\mathbb{R}^2$, complexify, then quotient by a naturally defined $\mathbb{Z}^2$-action.  For brevity, we describe here a single example, as full details appear in \cite{IW9}.

\begin{example}
Consider $E_8$, with choice of the following two shaded nodes:
\[
\begin{tikzpicture}[scale=0.8]
 \node (m1) at (-1,0) [DB] {};
 \node (0) at (0,0) [DW] {};
 \node (1) at (1,0) [DB] {};
 \node (1b) at (1,0.75) [DB] {};
 \node (2) at (2,0) [DB] {};
 \node (3) at (3,0) [DB] {};
 \node (4) at (4,0) [DB] {};
  \node (5) at (5,0) [DW] {};
\node (m1a) at (-1,-0.3) {$\scriptstyle 2$};
 \node (0a) at (0,-0.3) {$\scriptstyle 4$};
 \node (1a) at (1,-0.3) {$\scriptstyle 6$};
 \node (1ba) at (0.75,0.75) {$\scriptstyle 3$};
 \node (2a) at (2,-0.3) {$\scriptstyle 5$};
 \node (3a) at (3,-0.3) {$\scriptstyle 4$};
 \node (4a) at (4,-0.3) {$\scriptstyle 3$};
 \node (5a) at (5,-0.3) {$\scriptstyle 2$};
\draw [-] (m1) -- (0);
\draw [-] (0) -- (1);
\draw [-] (1) -- (2);
\draw [-] (2) -- (3);
\draw [-] (3) -- (4);
\draw [-] (4) -- (5);
\draw [-] (1) -- (1b);
\end{tikzpicture}
\]
A similar calculation, albeit slightly more combinatorially complicated than before, results in the following infinite real  hyperplane arrangement. 
\[
\begin{array}{c}
\includegraphics[angle=0,scale = 1]{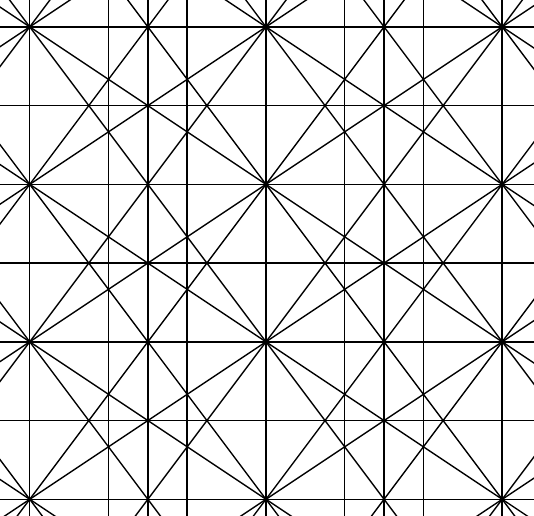}
\end{array}
\]
There is an obvious $\mathbb{Z}^2$-action.
\end{example}

\section{Relationship to Flops}

For simplicity in this section, consider only the case where $X\to \Spec R$ flops a single curve, namely the fibre with reduced scheme structure is
\[
\Curve\colonequals f^{-1}(\m)^{\rm red}\cong \mathbb{P}^1.
\]
The combinatorics in Section~\ref{L0} determine stability conditions on certain subcategories of $\Db(\coh X)$.  Details will appear in \cite{HW}. Here, we instead focus on the new autoequivalences that arise, and the resulting action of $\uppi_1(S^2\backslash\{N+2\})$ on $\Db(\coh X)$.  Full details of this will appear in \cite{DW5}. 

\subsection{General Elephants} As is well-known, considering the pullback via a generic hyperplane section $g\in R$ gives a commutative diagram
\[
\begin{tikzpicture}
\node (A1) at (0,0) {$Y$};
\node (A2) at (3,0) {$X$};
\node (B1) at (0,-2) {$\Spec (R/g)$};
\node (B2) at (3,-2) {$\Spec R$};
\draw[->] (B1) --(B2);
\draw[->] (A2) --node[right]{$\scriptstyle f$}(B2);
\draw[densely dotted,->] (A1) --(A2);
\draw[densely dotted,->] (A1) --node[left]{$\scriptstyle \upvarphi$}(B1);
\end{tikzpicture}
\]
where $R/g$ has only ADE surface singularities, and $\upvarphi\colon Y\to \Spec (R/g)$ is a \emph{partial} crepant resolution.  As such, $Y$ is dominated by the minimal resolution.

By McKay correspondence, the minimal resolution of $\Spec(R/g)$ is controlled by an ADE Dynkin diagram.  Since $Y$ is obtained from the minimal resolution by contracting curves, the combinatorics of $Y$ are controlled by a \emph{shaded} ADE Dynkin diagram, where we shade the vertices corresponding to the curves that are \emph{not} contracted in order to obtain $Y$.   Thus, the shaded vertices correspond to the curves that \emph{are} in $Y$, and so the number of shaded vertices equals the number of curves in $Y$, and thus in $X$.  By the assumption in this section, this number equals one.

\begin{defin}
The unique shaded vertex has a label $\ell$, in Section~\ref{S1.1}. This is called the \emph{length} of the curve.  
\end{defin}

It is a theorem of Katz--Morrison that when $X$ is smooth, only certain shaded arrangements can appear.  However, this fails in the singular setting described here. The combinatorial Remark~\ref{which ell not relevant} saves the day: because of it, length $\ell$ curves turn out to have the same stability conditions and autoequivalences, regardless in which Dynkin diagram they appear.

\subsection{Helices}
 Recall first that for a length $\ell$ flop, exactly as in \cite{Katz} there are sheaves
\[
\scrO_{\Curve},\scrO_{2\Curve},\hdots,\scrO_{\ell\Curve}
\]
on $Y$, and thus on $X$.    Each $a\Curve$ is a CM scheme of dimension one.  The dichotomy in part (2) of the following result will turn out to explain the jump in the combinatorics in Remark~\ref{which ell not relevant}.

\begin{prop}
For a length $\ell$ flop, the following hold.
\begin{enumerate}
\item $\omega_{2\Curve}\cong\scrO_{2\Curve}(-1)$.
\item $\Ext^1_X(\scrO_{2\Curve},\scrO_{3\Curve})=0\iff\ell\leq 4$.
 Otherwise $\Ext^1_X(\scrO_{2\Curve},\scrO_{3\Curve})=\mathbb{C}$.
 \end{enumerate}
\end{prop}
Thus, in the case $\ell=5,6$, there is a unique non-split extension
\[
0\to\scrO_{3\Curve}\to\scrZ\to\scrO_{2\Curve}\to 0. 
\]
Given this, we can now introduce our main new concept, that of a simples helix associated to a length $\ell$ flop.  Recall that, given $\ell$, there is an associated number $N$ in Remark~\ref{which ell not relevant}, which records the number of holes on the equator.
\begin{defin}
The \emph{simples helix} $\{\scrS_i\}_{i\in\mathbb{Z}}$ is defined to be the $\mathbb{Z}$-indexed family of sheaves, where for all $i$ 
\[
\scrS_{i+N}\cong\scrS_{i}\otimes\scrO(1),
\]
and $\scrS_0,\hdots,\scrS_{N-1}$ is defined to be
\[
\left\{
\begin{array}{ll}
\scrO_{\Curve}(-1)\,,\scrO_{\ell\Curve}\,,\hdots\,,\scrO_{3\Curve}\,,\scrO_{2\Curve}\,,\omega_{3\Curve}(1)\,,\hdots\,,\omega_{\ell\Curve}(1)&\mbox{if }\ell\leq 4\\
\scrO_{\Curve}(-1)\,,\scrO_{\ell\Curve}\,,\hdots\,,\scrO_{3\Curve}\,, \scrZ\,, \scrO_{2\Curve}\,,\scrZ^\omega(1)\,,\omega_{3\Curve}(1)\,,\hdots\,,\omega_{\ell\Curve}(1)&\mbox{if }\ell=5,6.
\end{array}
\right.
\]
Here $\scrZ^\omega$ is the sheaf such that $\mathbb{D}(\scrZ)\cong\scrZ^\omega[1]$, where $\mathbb{D}$ is duality.
\end{defin}

Consider $\scrA$, the category of perverse sheaves with perversity zero.  Suitably interpreted, this has two simples, which are known to be $\scrS_{-1}[1]$ and $\scrS_{0}$ \cite{VdB1d}.  We can tilt $\scrA$ at the second simple $\scrS_0$, and obtain a new heart $\scrA_1$.  We can then iterate: tilt $\scrA_1$ at its first simple to obtain $\scrA_2$, then tilt $\scrA_2$ at its second simple to obtain $\scrA_3$, etc.  This process gives a $\mathbb{Z}$-indexed family of hearts.

\begin{thm}
For a length $\ell$ flop, and for all $t\in\mathbb{Z}$, the following hold.
\begin{enumerate}
\item The simples of $\scrA_t$ are $\scrS_{t-1}[1]$ and $\scrS_t$.
\item $\scrA_t$ has a progenerator $\scrV_{t-1}\oplus\scrV_t$ which is a tilting vector bundle.
\item Noncommutative deformations of $\scrS_t$ are represented by an algebra $\Lambda_t^{\rm def}$, which is a factor of $\End_X(\scrV_t)$.
\end{enumerate}
\end{thm}

For our purposes later, the factor in part (3) is what allows us to control the noncommutative deformations: both from the viewpoint of autoequivalences in the next subsection, and from the viewpoint of curve-counting invariants in Section~\ref{applications}.

\subsection{Monodromy}
We claim that the simples helix in the previous subsection describes the monodromy of an action of $\uppi_1(S^2\backslash\{N+2\})$ on $\Db(\coh X)$.

Visually, $\uppi_1$ is generated as a group by $a,b_0,\hdots,b_{N-1},c$, where $a$ is monodromy around the north pole, $b_0,\hdots,b_{N-1}$ are mondromy around the holes on the equator, and $c$ is monodromy around the south pole.  Choosing orientations carefully (see \cite[Theorem 6.5]{DW5}), the relation
\[
c\circ b_0\circ\hdots\circ b_{N-1}\circ a=1
\]
holds.  Our main result does two things: it first constructs new autoequivalences in part (2), then shows that they describe monodromy in part (3).

\begin{thm}
For a length $\ell$ flop, and for all $t\in\mathbb{Z}$, the following hold.
\begin{enumerate}
\item Even although the sheaf $\scrS_t$ need not be perfect (as a complex), its universal sheaf $\scrE_t$ from noncommutative deformation theory is perfect.
\item There is a twist autoequivalence $\twistGen_{\scrS_t}$ of $\Db(\coh X)$ which fits into a functorial triangle
\[
\RHom_X(\scrE_t,-)\otimes_{\Lambda_t^{\rm def}}\scrE_t\to \Id\to\twistGen_{\scrS_t}\to
\]
\item For $i=0,\hdots, N-1$, the assignment
\[
\begin{array}{rcl}
a &\longmapsto&-\otimes\scrO(-1)\\
b_i&\longmapsto&\twistGen_{\scrS_i}\\
c&\longmapsto&\flop^{-1}\circ\big(-\otimes\scrO(-1)\big)\circ\flop,
\end{array}
\]
where $\flop$ is the flop functor, induces a group homomorphism 
\[
\uppi_1(S^2\backslash\{N+2\})\to\Auteq\Db(\coh X).
\] 
\end{enumerate}
\end{thm}
The general philosophy of String K\"ahler Moduli Spaces is that the action in part (3) should be faithful.

\subsection{Applications}\label{applications}
For brevity we describe only curve-counting applications.  More are given in \cite{DW5}.  Since curve invariants are only defined in the smooth setting, here we assume that $X$ is smooth.  Recall that the contraction algebra is defined to be $\Lambda_0^{\rm def}$, namely the algebra that represents noncommutative deformations of $\scrO_{\Curve}(-1)$.
\begin{thm}\label{main app}
Consider a length $\ell$ flop $X\to\Spec R$, where $X$ is smooth.  
\begin{enumerate}
\item For any $1\leq a\leq \ell$, the following conditions are equivalent.
\begin{enumerate}
\item Strictly noncommutative deformations of $\scrO_{a\Curve}$ exist. 
\item  $2a\leq \ell$.
\item Higher multiples of $a\Curve$ exist.
\end{enumerate}
\item
The lower bounds for the Gopakumar--Vafa invariants, and the lower bound for the dimension of the contraction algebra, are as follows.
\[
\begin{tabular}{clc}
\toprule
$\ell$&\textnormal{GV lower bound}&$\dim\CA$ \textnormal{lower bound}\\
\midrule
$1$ & $(1)$& $1$\\
$2$ & $(4,1)$& $8$\\
$3$ & $(5,3,1)$& $26$\\
$4$ & $(6,4,2,1)$& $56$\\
$5$ & $(7,6,4,2,1)$& $124$\\
$6$ &$(6,6,4,3,2,1)$& $200$\\
\bottomrule
\end{tabular}
\]

\end{enumerate}
\end{thm}

\begin{remark}
When $\ell=1,2,6$, it is known that the lower bound is realised.  At this stage, it remains unclear whether  the lower bound can be obtained when $\ell=3,4,5$. However, in those cases, examples of flops with GV invariants $(6,3,1)$, $(6,5,2,1)$ and $(8,6,4,2,1)$ are known to exist. 
\end{remark}

\end{document}